\documentclass[a4paper]{scrartcl}
\usepackage[]{amsmath,amssymb}
\usepackage{amsthm}
\usepackage{enumitem}
\usepackage[utf8]{inputenc}
\usepackage[T1]{fontenc}
\usepackage{lmodern}
\usepackage{fourier}
\usepackage{graphicx}
\usepackage{lmodern}
\usepackage{mathtools}
\usepackage{dsfont,mathrsfs}

\DeclareMathOperator{\oD}{d}
\def\D{\oD\!}
\def\dif{\,\D}

\usepackage{tikz}
\def\R{{\mathbf R}}

\def\N{{\mathbf N}}

\def\car#1{{\mathbf 1}}

\newcommand{\esp}[1]{{\mathbf E}\left[#1\right]}

\def\/{\,|\,} 

\def\car{\mathbf 1}

\usepackage{amsmath,amsfonts,amssymb,amsthm} 
\newtheoremstyle{hermesexercises}{11pt}{11pt}{\normalfont}{0pt}{\scshape}{.--}{.5em}{}
{\theoremstyle{hermesexercises}
}
\newtheoremstyle{hermesremark}{11pt}{11pt}{\normalfont}{0pt}{\scshape}{.--}{.5em}{}
\newtheoremstyle{myhermesremark}{11pt}{11pt}{\normalfont}{0pt}{\scshape}{.--}{.5em}{}
{\theoremstyle{myhermesremark}
  \newtheorem{remark}{Remark}

}
\newtheorem{notation}{Notation}
\newtheoremstyle{hermestheorem}{11pt}{11pt}{\normalfont}{0pt}{\scshape}{.--}{.5em}{}
{\theoremstyle{hermestheorem}
\newtheorem{theoremT}{Theorem}
\newtheorem{definitionT}{Definition}
\newtheorem{corollaryT}[theoremT]{Corollary}
\newtheorem{lemmaT}[theoremT]{Lemma}}
\usepackage{tcolorbox}
\tcbuselibrary{skins,xparse,breakable}
\definecolor{fonce}{HTML}{323031}
\definecolor{cadet}{HTML}{084C61}
\definecolor{vertdeau}{HTML}{177E89}
\definecolor{ocre}{HTML}{FFC857}
\definecolor{creme}{HTML}{DB3A34}

\tcbset{textmarker/.style={%
skin=enhancedmiddle jigsaw,breakable,parbox=false, boxrule=0mm,leftrule=5mm,rightrule=5mm,boxsep=0mm,arc=0mm,outer arc=0mm, left=3mm,right=3mm,top=1mm,bottom=1mm,toptitle=1mm,bottomtitle=1mm,oversize}}

\tcbset{textmarker2/.style={%
skin=enhancedmiddle jigsaw,breakable,parbox=false, boxrule=0mm,leftrule=3mm,boxsep=0mm,arc=0mm,outer arc=0mm, left=3mm,right=3mm,top=1mm,bottom=1mm,toptitle=1mm,bottomtitle=1mm,oversize}}

\tcolorboxenvironment{lemmaT}{textmarker2,colback=ocre!5!white,colframe=ocre}
\newenvironment{lemma}{%
     \begin{lemmaT}}{\end{lemmaT}}

   \tcolorboxenvironment{definitionT}{textmarker2,colback=creme!5!white,colframe=creme}
   \newenvironment{definition}{%
     \begin{definitionT}}{\end{definitionT}}
      \tcolorboxenvironment{theoremT}{textmarker2,colback=vertdeau!5!white,colframe=vertdeau}
 \newenvironment{theorem}{%
     \begin{theoremT}}{\end{theoremT}}
      \tcolorboxenvironment{corollaryT}{textmarker2,colback=vertdeau!5!white,colframe=vertdeau}
 \newenvironment{proposition}{%
     \begin{theoremT}}{\end{theoremT}}
      \tcolorboxenvironment{corollaryT}{textmarker2,colback=vertdeau!5!white,colframe=vertdeau}
      \newenvironment{corollary}{%
     \begin{corollaryT}}{\end{corollaryT}}
\usepackage{enumitem}
\usetikzlibrary{%
  arrows,%
  calc,%
  fit,%
  patterns,%
  plotmarks,%
  shapes.geometric,%
  shapes.misc,%
  shapes.symbols,%
  shapes.arrows,%
  shapes.callouts,%
  shapes.multipart,%
  shapes.gates.logic.US,%
  shapes.gates.logic.IEC,%
  er,%
  automata,%
  backgrounds,%
  chains,%
  topaths,%
  trees,%
  petri,%
  mindmap,%
  matrix,%
  calendar,%
  folding,%
  fadings,%
  through,%
  positioning,%
  scopes,%
  decorations.fractals,%
  decorations.shapes,%
  decorations.text,%
  decorations.pathmorphing,%
  decorations.pathreplacing,%
  decorations.footprints,%
  decorations.markings,%
  shadows}
\newcommand{\iid}{\textnormal{i.i.d.\ }}

\renewcommand{\esp}{\mathbb{E}}

\newcommand{\prob}{\mathbb{P}}
\newcommand{\ind}{\mathds{1}}
\newcommand{\bs}{\setminus}
\newcommand{\lb}{\{}
\newcommand{\rb}{\}}

\newcommand\numberthis{\addtocounter{equation}{1}\tag{\theequation}}

\newcommand{\dom}{\operatorname{Dom}}
\newcommand{\lipbis}{\mathrm{Lip}_{{[2]}}(\R)}
\def\eqdis{\overset{\mathrm{d}}{=}}

\def\Po{\mathbf{P}^{0}}
\def\Lo{\mathscr{L}_{0}}

\def\myesp#1{\mathbf{E}\left[#1\right]}
\def\Pa{\mathbf{P}^{\alpha}} 

\title{Convergence rate for the coupon collector's problem with Stein's method}
\author{B. Costacèque \and L. Decreusefond}
\date{}
\begin{document}
\maketitle{}
\begin{abstract}
    The functional characterization of a measure, an essential but delicate
    aspect of Stein's method, is shown to be accessible for stable probability
    distributions on convex cones. This notion encompasses the usual stable
    distributions \textit{e.g.} Gaussian, Pareto, \textit{etc.} but also the
    max-stable distributions: Weibull, Gumbel and Fréchet. We use the definition
    of max-stability to define a Markov process whose invariant measure is the
    stable measure of interest. In this paper, we focus on the Gumbel
    distribution and show how this construction can be applied to estimate the
    rate of convergence in the classical coupon collector's problem.
\end{abstract}

\noindent{Math Subject Classification: 60F05,60E07}
\def\LL{\operatorname{L}}
\def\laurent#1{\footnote{\textcolor{blue}{#1}}}
\section{Introduction}
Stein's method in its modern acceptance (see \cite{Decreusefond15,
  Decreusefond22}) is based upon an idea the physicists call \textit{stochastic
  quantization}. Given a probability measure~$\mu$ on a space $(E,{\cal E})$,
this procedure established in \cite{Parisi1981} consists in constructing a
Markov semi-group $(P_{t})_{t\geq 0}$ whose invariant distribution is~$\mu$. In
physics, this approach is, for instance, used for numerical simulations of gauge
theories with fermions \cite{namiki1992stochastic}. In the context of the
Stein's method, we use this construction to establish the representation
formula:
\begin{equation}
  \label{eq_dalphas_core:1}
  \int_{E} f\dif \mu-\int_{E} f\dif \nu=\int_{E} \int_{0}^{\infty}LP_{t}f \dif t
  \dif \nu,
\end{equation}
where $L$ is the generator of the semi-group and $\nu$ is any other probability
measure on $(E,{\cal E})$. Alternatively, $L$ can be defined by the property:
\begin{equation}
  \label{eq_dalphas_core:3}\int_{E}Lf\dif \nu =0, \ \forall f \in \dom L \Longleftrightarrow \nu=\mu
\end{equation}
Describing $L$ or $(P_{t})_{t\ge 0}$ is equivalent in the sense that we do have
formulas which link one to the other:
\begin{equation*}
  \left.Lf=\frac{\dif}{\dif t}P_{t}f\right|_{t=0}  \text{ and } \frac{\dif}{\dif t}P_{t}f=LP_{t}f.
\end{equation*}
Actually, the computations may not be so straightforward.
So, shall we start with the operator $L$ or with the Markov
semi-group?
There is no magical recipe to construct them for an arbitrary probability
measure. See~\cite{Ley_2017} and references therein for a thorough survey of how
to construct an operator satisfying~\eqref{eq_dalphas_core:3} for a wide range
of probability measures. Note that the operator so obtained may not be the
generator of a Markov process. For the standard Gaussian measure on $\R^{n}$,
the usual operator is
\begin{equation*}
  Lf(x)=-\langle x,\ \nabla f(x)\rangle + \Delta f(x),
\end{equation*}
where $\nabla f(x) = (\partial_{1}f(x),\dots, \partial_{n}f(x))$ and
$\Delta f(x) = \sum_{j=1}^{n}\partial^{2}_{j,j}f(x)$. The corresponding
semi-group is the so-called Ornstein-Uhlenbeck semi-group is given by
\begin{equation*}
  P_{t}f(x)=\int_{\R^{n}}f\Bigl(e^{-t}x+\sqrt{1-e^{-2t}}y\Bigr)\dif \mu(y).
\end{equation*}
If $\mu$ denotes the law of a Poisson process $N$ of control measure $\sigma$,
the so-called Glauber semi-group is defined by
\begin{equation*}
  P_{t}f(\omega)=\mathbb{E}\Big[f\big(e^{-t}\circ \omega+(1-e^{-t})\circ N\big)\Big]
\end{equation*}
where $p\circ \omega$ is the $p$-thinning of the configuration $\omega$ (see
\cite{Decreusefond16}). The crucial remark is that the semi-group property of
the operators defined above is a mere rewriting of the stability of the Gaussian
and Poisson measures: The standard Gaussian measure on $\R^{n}$ is the unique
probability law on $\R^{n}$ such that
\begin{equation}
  \label{eq_coupon_core:3}
  e^{-t}\,Z'+\sqrt{1-e^{-2t}}\,Z''\eqdis Z
\end{equation}
where $Z',Z''$ are two independent copies of $Z$. The Poisson process of control
measure $\sigma$ is characterized by the fact that
\begin{equation}
  \label{eq_dalphas_core:4}e^{-t}\circ N' + (1-e^{-t})\circ N''\eqdis N
\end{equation}
where $N',N''$ are two independent copies of $N.$ This invariance also entails
that the Gaussian (respectively the Poisson) measure is the stationary measure
of the semi-group $(P_{t})_{t\ge 0}$. Remark that at a formal level,
equations~\eqref{eq_coupon_core:3} and~\eqref{eq_dalphas_core:4} can be written
as
\begin{equation}\label{stability}
  D_{a}Z'\oplus D_{b}Z''\eqdis Z
\end{equation}
where $a^{2}+b^{2}=1$ in~\eqref{eq_coupon_core:3}, $a+b=1$
in~\eqref{eq_dalphas_core:4} and $\eqdis$ denotes an equality in distribution.
The operator $D_{a}$ is the multiplication (respectively the thinning) of real
numbers (respectively of point processes) and the sign $\oplus$ represents the
addition of real numbers (respectively the superposition of point processes)
in~\eqref{eq_coupon_core:3} (respectively~\eqref{eq_dalphas_core:4}). This means
that the algebraic structure at play here is that of a semi-group (the plus
operation) with a family $(D_{a})_{a\in T}$ of automorphims which satisfy the
identity $D_{a}\circ D_{b}=D_{ab}$; the transformation $D_{a}$ is called a
\textit{multiplication} (by $a$). In the seminal paper~\cite{Davydov2008}, such
a structure is called a \textit{convex cone}. It is shown there that there exist
numerous other examples of stable distributions provided that we change the
meaning of the $\oplus$ sign and the group of automorphisms. These distributions
are often interesting because their stability implies that they appear in limit
theorems akin to the central limit theorem.

In this paper, we focus on the Gumbel distribution which is characterized by the
stability identity
\begin{equation}
  \label{max_stab}
  \max\Bigl(Z'+\log a,\ Z''+\log(1-a)\Bigr) \eqdis Z,
\end{equation}
for any $a\in (0,1)$. This falls in the previously introduced framework where
$x\oplus y=\max(x,y)$ and $D_{a}x = x+\log a$. We then take advantage of this
identity to define a Markov semi-group which leaves the Gumbel distribution
invariant. We then compute its generator $L$ and show that the
equivalence~\eqref{eq_dalphas_core:3} holds.

Let $I_{n} \coloneq \{1,\dots,n\}$ and consider $(X_{k})_{k \geq 1}$ a sequence
of independent random variables, uniformly distributed over~$I_{n}$. We denote
by
\begin{equation*}
  T_{n} \coloneq \inf\big\{k\ge 1, \{X_{1},\dots,X_{k}\}=I_{n}\big\}.
\end{equation*}
The analysis of the properties of $T_{n}$ is the so-called coupon collector's
problem. It is well known that
\begin{equation}
  \label{eq_coupon_core:5}Z_{n}\coloneq \frac{T_{n}}{n}-\log n \xrightarrow{\text{law}}\mathcal{G}(0,1),
\end{equation}
where $\mathcal{G}(0,1)$ denotes the standard Gumbel distribution, with
cumulative distribution function $\exp(-\exp(-x))$ on $\R$.

In this paper, we use our semi-group and generator, which are characteristic to
the Gumbel distribution, to develop Stein's method in order to assess the order
of convergence in~\eqref{eq_coupon_core:5} for a Wasserstein-like distance. The
proof of the bound of the convergence rate is unusually lengthy. However, we
believe that it is inherently difficult to compress further, as it introduces novel techniques of
independent interest. The primary challenge arises from the fact that, unlike
the classical central limit theorem (CLT), we do not work with the partial sum of a sequence of
random variables but rather with a triangular array of independent but not
identically distributed random variables whose contributions vary in order of
magnitude. Consequently, the leave-one-out method and standard Taylor expansions
—commonly employed in such contexts— are not directly applicable here.
Furthermore, the generator associated with the Gumbel distribution presents
additional computational difficulties, as its form is less tractable than in
other settings. To address this, we develop several preliminary lemmas that
transform the term $\myesp{f'(Z_{n})}$ into a form amenable to comparison with
the remaining components of the generator.

To the best of our knowledge, very few papers so far have given uniform bounds
for the coupon collector problem: \cite{Sandor93} seems to be the first and the
convergence rate is the same as ours but for the Kolmogorov distance and through
a different method, which is by no means elementary. These results have been
extended to generalizations (where one considers the time $T_{n,m}$ needed to
$n-m$ distinct coupons for the first time, with $0 \le m \le n-1$) of the
original coupons collector problem in~\cite{Posfai10}. In particular, it is proved
therein that, as $n$ goes to infinity,
\[
  \underset{x \in \R}{\sup} \Big\vert \prob\Big(\frac{T_{n}}{n}-\log n \leq x\Big) - e^{-e^{-x}} \Big\vert = O\Big(\frac{\log n}{n}\Big),
\]
implying that the speed of convergence $n^{-1}\log n$ is optimal, at least in
Kolmogorov distance. A comprehensive and more recent reference regarding
non-uniform asymptotic expansions in coupon collector problems is \cite[Section
4.3]{Mladenovic24}.

Besides the bound for the original coupon collector problem, another
contribution of this paper is the introduction of the semi-group
$(\Po_{t})_{t\geq 0}$, making possible to use the generator approach of Stein's
method when the target distribution is the Gumbel measure. Let us mention that
Stein's method has been previously applied to extreme value distributions:
through point process techniques and Stein's method for the Poisson process,
\cite{Cipriani16,Feidt10} managed to derive bounds for the distance between a
maximum of (not necessarily independent) discrete random variables and extreme
value distributions. However their method did not rely on a characterizing
operator tailored for the target max-stable distribution. On that regard,
\cite{Bartholome13} is the first to define a Stein operator for the Fréchet
distribution and use it to find a bound for the extreme value theorem. Similar
characterizing operators have been found for the Gumbel and Weibull
distributions in \cite{Kusumoto20}. One drawback of the approach of
\cite{Bartholome13} is that one needs to assume the approaching random variable
to have a density, so as to work with its score function, thus excluding the
possibility of dealing with discrete random variables as we do in this paper.
Moreover, their approach to find a characterizing operator relies on the density
of the Fréchet distribution and so does not generalize well to max-stable random
vectors, as they do not always admit a density. We explain below how our
construction extends to the multivariate setting. To the best of our knowledge,
those constitute the main contributions to the problem of finding uniform rate
of convergence to extreme value distributions by using Stein's method.

It is possible to extend our results in several directions: first by considering
one of the many generalizations of the coupon collector problem which involve
extreme value distributions. For instance, one could study the typical time
$T_{n,m}$ needed to obtain $n-m$ distinct items for the first time. Another
classic variation of this problem asks for the typical behavior of the waiting
time $T_{n}^{(r)}$ for all coupons to appear $r$ items at least, for some
integer $r$. The special case $r = 1$ corresponds to the classical coupon
collector problem. For an arbitrary $r$, up to an appropriate normalization,
the waiting time converges to a Gumbel distribution as well, making it possible
to apply the tools we introduced to it. The first variation should be easily
accessible to our method, as both $T_{n}$ and $T_{n,m}$ can be expressed as sums
of independent random variables. As far as we know, this is not true anymore for
$T_{n}^{(r)}$, and thus a different approach is required.
  
Second, one could try and work with a different probability metric
$\operatorname{d}$ than the smooth Wasserstein distance we use. This raises
several difficulties however, as the test-functions $h$ associated with
$\operatorname{d}$ must be regular enough so that $\Po_{t}h$ is differentiable.
For instance, if we assume $h$ to be just $1$-Lipschitz, then $\Po_{t}h$
may not be differentiable everywhere. This is problematic when working with
discrete random variables, such as $T_{n}$. That issue is common when applying
Stein's method, and several fixes have been developed, such as bounding the
target metric by a function of another, more amenable metric. This solution is
generic although the resulting bound will be unlikely to be optimal. We refer to
\cite{Gaunt23} for a thorough presentation.

Lastly, one could change the definition of $\Po_{t}$ to define a Markov
semi-group for the two remaining extreme value distributions: the Weibull
distribution and the Fréchet distribution. Indeed they both satisfy a stability
relation akin to \eqref{stability}. For instance, if $\alpha > 0$, the Fréchet
distribution $\mathcal{F}(\alpha)$ is such that:
\[
  aZ' \oplus bZ'' \eqdis (a^{\alpha} + b^{\alpha})Z,\ a,b \geq 0
\]
where $Z'$ and $Z''$ are \iid copies of $Z \sim \mathcal{F}(\alpha)$ and
$x \oplus y = \max(x,y)$. One can then define a family linear operators
$(\Pa_{t})_{t \geq 0}$ by setting
\[
  \Pa_{t}f(x) \coloneq \esp\Big[f\big(e^{-\frac{t}{\alpha}}x \oplus (1-e^{-t})^{\frac{1}{\alpha}}Z\big)\Big]
\]
and prove that this is indeed a semi-group admitting $\mathcal{F}(\alpha)$ as a
stationary measure. This idea readily extends to the Weibull distribution, and
even to higher-dimensional extreme value distributions, the so-called max-stable
random vectors. We refer the reader to \cite{Costaceque24_mlti, Costaceque24}
for more about those generalizations and their applications to Stein's method.

This work is organized as follows. In the second Section, we introduce a new
Markov semi-group $(\Po_{t})_{t\geq 0}$ whose stationary measure is the Gumbel
distribution and give its essential properties. In particular, we obtain a Stein
identity characterizing the Gumbel distribution. The third Section contains the
main result of this paper, that is, a rate of convergence for the standard
coupon collector problem in smooth Wasserstein distance. The last Section is
dedicated to the proof of that result.

\def\lip{\operatorname{Lip}} \def\lipbis{\mathcal{C}^{1,\lip}}
\section{Stochastic quantization of the Gumbel distribution}

The Gumbel distribution is one of the three so-called \textit{extreme value
  distributions}. It arises in the study of extreme events, for instance when
one aims at predicting record floods or drought. The most famous occurrence of
the Gumbel distribution is given by the celebrated Fisher-Tippett-Gnedenko
theorem \cite{Resnick87}, which states that a renormalized maximum of
independent and identically distributed (iid) random variables can have only
three possible limits: the Fréchet distribution, the (negative) Weibull
distribution and the Gumbel distribution. Which one appears depends on the right
tail of the initial distribution. The Fréchet distribution corresponds to
regularly varying tails, the Weibull distribution to bounded right tails, and
the Gumbel distribution to light tails. The domain of attraction of the latter
thus includes many fundamental distributions, like the Gaussian distribution,
the log-normal distribution, the exponential distribution, the logistic
distribution, {etc.}

The set of Lipschitz continuous from $\R$ to $\R$ is denoted by $\lip$ and is
equipped with the norm
\begin{equation*}
  \|f\|_{\lip}=\sup_{t\neq s} \frac{|f(t)-f(s)|}{|t-s|}\cdotp
\end{equation*}
We denote by $\lipbis $, the class of Lipschitz continuous functions whose first
derivative is also Lipschitz over $\R$:
\begin{equation*}
  \lipbis \coloneq \big\lb f \in \lip,\  f'\text{ admits a representative in } \lip\big\rb.
\end{equation*}
We equip $\lipbis$ with the semi-norm
\begin{equation*}
  \|f\|_{\lipbis}=\max(\|f\|_{\lip},\, \|f'\|_{\lip}).
\end{equation*}
Note that we do not impose $f\in \lipbis$ to be bounded.
\begin{definition}
  Let $\mu$ be the Gumbel distribution, with p.d.f. $x \mapsto e^{-(x+e^{-x})}$
  on $\R$ and $Z$ a random variable with that distribution. Let
  $f \in \lipbis $. Recall that $x \oplus y = \max(x,y)$ and set for all
  $x \in \R$:
	$$
	D_{a}x \coloneq x + \log a,\ a \geq 0
	$$
	with $D_{0}x = -\infty$, and define a family of linear operators on
    $\lipbis $ by setting:
	\begin{align*}
	  \Po_{t}f(x) &\coloneq \esp\big[f\big(D_{e^{-t}}x \oplus D_{1-e^{-t}}Z\big)\big]\\
				  &=\esp\left[f\Bigl(\max\bigl(x-t,\ Z + \log(1-e^{-t})\bigr)\Bigr)\right],
	\end{align*}
	for all $t \in \R_{+}$ and $x \in \R$.
  \end{definition}
  \begin{proposition}
    \label{thm_coupon_core:characterization}
	The family of operators $\Po=(\Po_{t})_{t\geq 0}$ is an ergodic semi-group.
    It is stationary for the Gumbel distribution $\mu$.
  \end{proposition}
  \begin{proof}

	It is clear that $\Po_{0} = \operatorname{Id}$ and $\Po_{t}$ is a linear
    operator well defined on $\lipbis
    $ 
	for all $t\ge 0$.

	Denote by $Z,\, Z'$ two independent random variables with Gumbel
    distribution. For all $t,s \geq 0$:
	\begin{multline*}
	  (\Po_{t}\circ \Po_{s})f(x) = \Po_{t}(\Po_{s}f)(x)\\
	  \begin{aligned}
		&= \esp\Big[\Po_{t}f\big((x - s)\oplus (Z + \log(1-e^{-s}))\big)\Big]\\
		&= \esp\Big[f\big((x - s - t)\oplus (Z + \log(1-e^{-s}) - t) \oplus (Z' + \log(1-e^{-t}))\big)\Big]\\
		&= \esp\Big[f\big((x - s - t)\oplus (Z + \log(e^{-t}(1-e^{-s}))) \oplus (Z' + \log(1-e^{-t}))\big)\Big]\\
		&= \esp\Big[f\big((x - s - t)\oplus (Z + \log(1-e^{-s-t}))\big)\Big]\\
		&= \Po_{t+s}f(x),
	  \end{aligned}
	\end{multline*}
	where the penultimate equality comes from~\eqref{max_stab}. From the
    definition of the maximum, we know that:
	\begin{multline}
	  \label{MaxDec}
	  \Po_{t}f(x) = f(x-t)\prob\big(Z + \log(1 - e^{-t})\leq x - t\big)\\
	  \shoveright{+ \esp\big[f\big(Z + \log(1-e^{-t})\big)\ind_{\lb Z + \log(1-e^{-t}) \geq x - t \rb}\big]}\\
	  = f(x-t)e^{-\gamma_{t}e^{-x}} + \esp\big[f\big(Z + \log(1-e^{-t})\big)\ind_{\lb Z \geq x - \log\gamma_{t} \rb}\big],
	\end{multline}
	where $\gamma_{t} \coloneq e^{t}-1$.
	Since $f$ is Lipschitz, we see by dominated convergence that
	\begin{equation*}
	  \forall x \in \R,\ \lim_{t\to \infty}\Po_{t}f(x)=\int_{\R}f\dif \mu,
	\end{equation*}
	hence the ergodicity of $(\Po_{t})_{t\geq 0}$. Finally, we have:
	\begin{align*}
	  \esp\big[\Po_{t}f(Z)\big] &= \esp\Big[f\big((Z-t)\oplus (Z + \log(1-e^{-t}))\big)\Big]\\
								&= \esp\Big[f\big(Z + \log(e^{-t} + 1 - e^{-t})\big)\Big] \\
								&= \esp[f(Z)],
	\end{align*}
	which means that $\mu$ is the invariant measure of the semi-group $\Po.$
  \end{proof}
  \begin{proposition}\label{sec2_gnrtr}
	The generator $\Lo $ of the semi-group $(\Po_{t})_{t\geq 0}$ is given for
    all $f \in \lipbis$ by:
	\begin{align}
      \Lo f(x)& = -f'(x) + e^{-x}\esp\big[f(x + Y) - f(x)\big] \label{eq_coupon_core:7}\\
			  &= -f'(x) + e^{-x}\esp\big[f'(X + Y)\big],\label{eq_coupon_core:8}
    \end{align}
	where $Y$ is a random variable with exponential distribution of
    parameter~$1$, denoted by $\mathcal{E}(1)$.
  \end{proposition}
  \begin{proof}
    First, we rewrite \eqref{MaxDec} as:
    \begin{multline}\label{eq_coupon_core:1}
      \Po_{t}f(x)  = f(x-t)e^{-\gamma_{t}e^{-x}}
      { + \int_{x - \log\gamma_{t}}^{\infty}f\big(z + \log\gamma_{t} - t\big)e^{-(z + e^{-z})}\dif z }\\
      = f(x-t)e^{-\gamma_{t}e^{-x}} + \gamma_{t}\int_{x}^{\infty}f(z - t)e^{-(z + \gamma_{t}e^{-z})}\dif z .
    \end{multline}
    A dominated convergence argument and differentiating with respect to $t$
    yield:
    \begin{multline}
      \label{sec2_3parts}
      \frac{\dif}{\dif t}\Po_{t}f(x)  = -f'(x-t)e^{-\gamma_{t}e^{-x}}- e^{t}f(x-t)e^{-x}e^{-\gamma_{t}e^{-x}} \\
      \shoveright{ + e^{t}\int_{x}^{\infty}f(z - t)e^{-(z + \gamma_{t}e^{-z})}\dif z+ \gamma_{t}\int_{x}^{\infty}\frac{\dif }{\dif t}\big(f(z - t)e^{-(z + \gamma_{t}e^{-z})}\big)\dif z}\\
      \shoveleft{  = -f'(x-t)e^{-\gamma_{t}e^{-x}} - e^{t}f(x-t)e^{-x}e^{-\gamma_{t}e^{-x}}}\\
      + e^{t}\int_{x}^{\infty}f(z - t)e^{-(z + \gamma_{t}e^{-z})}\dif z+  R_{t},
    \end{multline}
    where
    \begin{align*}
      R_{t} &\coloneq -\gamma_{t}\Big( \int_{x}^{\infty}f'(z - t)e^{-(z + \gamma_{t}e^{-z})}\dif z + e^{t}\int_{x}^{\infty}f(z - t)e^{-(2z + \gamma_{t}e^{-z})}\dif z \Big).
    \end{align*}
    That remainder is bounded by:
    \begin{align*}
      \vert R_{t} \vert &\leq \gamma_{t}\int_{x}^{\infty}e^{-(z + \gamma_{t}e^{-z})}\dif z + \gamma_{t}e^{t}\int_{x}^{\infty}\big(\vert f(0) \vert + \vert z \vert + t \big)e^{-(2z + \gamma_{t}e^{-z})}\dif z\\
                        &\leq \gamma_{t}\int_{x}^{\infty}e^{-z}\dif z + \gamma_{t}e^{t}\int_{x}^{\infty}\big(\vert f(0) \vert + \vert z \vert + t \big)e^{-2z}\dif z.
    \end{align*}
    Each of those integrals is finite for $x \in \R$ and $\gamma_{t}$ vanishes
    as $t$ goes to $0$, hence
    \begin{equation*}
      R_{t}\xrightarrow[t\to 0]{} 0.
    \end{equation*}
    On the other hand, it is easily shown that the first three terms in
    \eqref{sec2_3parts} have a non-trivial limit, given by:
    \begin{align}
      \Lo f(x) & = -f'(x) - e^{-x}f(x) + \int_{x}^{\infty}e^{-z}f(z)\dif z             \label{eq_coupon_core:9}   \\
			 & = -f'(x) + e^{-x}\int_{0}^{\infty}e^{-y}\big(f(x + y) - f(x)\big)\dif y,\notag
\end{align}
  hence \eqref{eq_coupon_core:7} holds.
  Furthermore, an integration by parts in~\eqref{eq_coupon_core:9} yields:
  \begin{align*}
	\Lo f(x) & = -f'(x) + \int_{x}^{\infty}e^{-y}f'(y)\dif y                                                   \\
			 & = -f'(x) + e^{-x}\int_{0}^{\infty}e^{-y}f'(x+y)\dif y.
  \end{align*}
  The last identity can be rewritten as in~\eqref{eq_coupon_core:8}.
\end{proof}
We now show that $\Lo$ yields a functional characterization of the Gumbel distribution.
\begin{theorem}[Stein's Lemma]
  Let $X$ be a random variable on $\R$ such that $\esp[e^{-X}]$ is finite. Then
  $X$ has the Gumbel distribution $\mathcal{G}(0,1)$ if and only if
  \begin{align}\label{cpn_Stein_lemma}
	\esp[f'(X)] = \esp\big[e^{-X}f'(X+Y)\big],
  \end{align}
  for all $f \in \lipbis$, where $Y$ has the exponential distribution
  $\mathcal{E}(1)$ and is independent of $X$.
\end{theorem}

\begin{proof}
   If $X$ has the Gumbel distribution, the result is an immediate consequence of Theorem~\ref{thm_coupon_core:characterization}.

  Conversely, assume that~\eqref{cpn_Stein_lemma} holds.  We first prove that the Laplace
  transform $L_{X}$ of $X$ is finite on $(-1,+\infty)$.
  By assumption, we already know that $L_{X}(1)$ is
  finite, and thus $L_{X}$ is well-defined on $[0,1]$. Now
  apply~\eqref{cpn_Stein_lemma} to
  $(x \mapsto e^{-tx})$, for some $t\in (-1,0]$, to obtain:
  \begin{equation}
	\label{eq_coupon_core:10}
	L_{X}(t) 
	= \frac{1}{1+t}L_{X}(1+t)\cdotp
  \end{equation}
	Hence, $L_{X}$ is finite on $(-1,1]$.  By induction,
	\eqref{eq_coupon_core:10} also  shows
	that $L_{X}$ is finite on $(-1,\infty)$ and
	 implies that
	\begin{equation*}
	   L_{X}(\lambda)=\lambda L_{X}(\lambda-1)
	 \end{equation*}
	 for all $\lambda>0$. Apply the Bohr-Mollerup Theorem  (see
	 \cite{Artin64}) to $\lambda \mapsto L_{X}(\lambda-1)$ to prove that
	 $L_{X}(\lambda)=\Gamma(\lambda+1)$, which actually means that $L_{X}$ is
	 equal  to the
	 Laplace transform of the Gumbel distribution.
  \end{proof}
  Next, we wish to determine how much $\Po_{t}$ leaves $\lipbis$ stable.
  \begin{corollary}\label{SeGr}
	Let $f \in \lipbis $. Then for any $t\ge 0$, $\Po_{t}f$ belongs to
	$\lipbis$: $\Po_{t}f$ is
	1-Lipschitz and its derivative is $2$-Lipschitz. Furthermore $\Lo \Po_{t}F$
	is well defined and for $t>0$, we have:
	$$
	\Vert \Lo \Po_{t}f \Vert_{\infty} \leq  \frac{2}{1-e^{-t}}
	\Vert f' \Vert_{\infty}.
	$$
  \end{corollary}
  \begin{proof}
It is an immediate consequence of~\eqref{eq_coupon_core:1} that for all
$x\in \R$ and all $t\ge 0$,
	\begin{align}\label{P_tDiff}
	  (\Po_{t}f)'(x) = f'(x-t)e^{-\gamma_{t}e^{-x}}.
	\end{align}
The first assertion follows  from~\eqref{P_tDiff}. The second claim comes from~\eqref{eq_coupon_core:8}:
   \begin{align*}
   	\int_{x}^{\infty}e^{-z}\vert(\Po_{t}f)'(z)\vert \dif z & \leq \Vert f' \Vert_{\infty}\int_{x}^{\infty}e^{-z}e^{-\gamma_{t}e^{-z}}\dif z \\
   														   & \leq \Vert f' \Vert_{\infty}\int_{\R}e^{-z}e^{-\gamma_{t}e^{-z}}\dif z \\
   														   & = \frac{1}{\gamma_{t}}\Vert f' \Vert_{\infty}
   \end{align*}
   for any $x \in \R$. This shows that $\Lo \Po_{t}f$ is defined, and even bounded as soon as $t$ is
   positive.
\end{proof}

\section{Rate of convergence to the Gumbel distribution}
Recall the notations introduced above. We are given a collection of $n$~images
and we draw one of each at random at each step. We denote by $T_{n}$ the number
of drawings before the completion of the collection. It is well-known that we
can write $T_{n}$ as
\begin{equation*}
  T_{n}=\sum_{i=1}^{n}\tau_{i}^{n}
\end{equation*}
where $(\tau_{i}^{n})_{1 \leq i \leq n}$ is a sequence of independent geometric
random variables of respective parameter $(n-i+1)/n$ with support on $\N^{*} = \lb 1,2, \dots \rb$. Furthermore, we set
\begin{equation*}
  Z_{n}=\frac{T_{n}-n\log n}{n}
\end{equation*}
and we want to estimate
\begin{equation*}
\operatorname{dist}_{\lipbis}\left(\operatorname{law}(Z_{n}),\ \mu\right)\coloneq\sup_{f\, :\, \|f\|_{\lipbis}=1}\left(\esp[f(Z_{n})]-\esp[f(Z)]\right).
\end{equation*}
Note that since $-f$ belongs to $\lipbis$ if $f \in \lipbis$, there is no need for absolute values over $\esp[f(Z_{n})]-\esp[f(Z)]$ in that definition. As  in	\cite[Theorem 4.1]{Coutin2013}, we can prove that the topology on the set of
probability measures on $\R$ induced by that distance is the same as the usual
Wasserstein-1 (or Kantorovitch-Rubinstein) topology.

The main theorem of this section is the following:
\begin{theorem}\label{cpn_main}
  Let $Z$ be a random variable with the standard Gumbel distribution
  $\mathcal{G}(0,1)$. Then there exists a constant $C>0$ such that for  $n \geq 2$, the following inequality holds:
  \begin{equation*}
\operatorname{dist}_{\lipbis}(\operatorname{law}(Z_{n}),\ \operatorname{law}(Z)) \leq C\,\frac{\log n}{n}\cdotp
  \end{equation*}
  \end{theorem}
The proof is lengthy and technical.
Before going further, we need to introduce some notations:
\begin{notation}
  Set
\begin{equation*}
  \alpha_{n} \coloneq 1-\frac{1}{n},\ \beta_{n} \coloneq -(n-1)\log \alpha_{n}, \ K_{n}\coloneq n\Big(1-\frac{1}{n}\Big)^{(n-1)\log(n-1)},
\end{equation*}
 as well as:
$$
\delta_{n} \coloneq \log\Big(1-\frac{1}{n}\Big) - \frac{1}{n}\log(n-1).
$$
We will make constant use of the following inequalities in the sequel for
$n\geq 2$:
\begin{align}\label{cpn_ineq_cstt}
  -1 \leq \delta_{n} \leq -\frac{1}{n} \leq 0 \leq \alpha_{n} < \beta_{n} < 1 \le  K_{n}.
\end{align}
\end{notation}
Using the results exposed in the previous part, we can initiate the generator
approach of the Stein's method. We know from \cite{Decreusefond15} that we can
write  for all $f \in \lipbis$
\begin{equation*}
\esp[f(Z_{n})] - \esp[f(Z)]  = -\int_{0}^{\infty}\esp\big[\Lo \Po_{t}f(Z_{n})\big] \dif t.
\end{equation*}
The unusual form of $\Lo$ (see~\eqref{eq_coupon_core:8}) makes non operative the
standard methods like leave-one-out or more generally integration by parts as in~\cite{Decreusefond19}.
 The main step of the proof is to establish the following
 identity~\eqref{cpn_IPP} which is very close to what is expected in view of~\eqref{cpn_Stein_lemma}. In the sequel $\operatorname{Geom}(p)$ will denote the positive geometric distribution, where $p \in (0,1]$.
\begin{lemma}
\label{thm_coupon_core:main_identity}
  Let $f$ be a function in $\mathscr{C}^{1}(\R, \R)$ and $n \geq 2$. We have the
  identity:
  \begin{align}\label{cpn_IPP}
	\esp\big[f'(Z_{n})\big] & = K_{n}\ \esp\big[e^{-\beta_{n}Z_{n-1}}f'\big(\alpha_{n}Z_{n-1} + {G_{n}} + \delta_{n}\big)\big],
  \end{align}
  where the random variable $G_{n}$ has distribution
  $n^{-1}\operatorname{Geom}(\frac{1}{n})$ and is independent of $Z_{n-1}$.
\end{lemma}
As a result, we see that we will have to compare $Z_{n}$ and $Z_{n-1}$ at several occasions.
 This is done by a very natural coupling.
\begin{lemma}\label{cpn_cplng}
  For every $n \in \N^{*}$, let $(\tau_{i}^{n})_{1 \le i \le n}$ be a sequence of positive \iid random variables with geometric distribution $\operatorname{Geom}(p_{i}^{n})$, with
  	$$
	p_{i}^{n} \coloneq \frac{n-i+1}{n},\ n \geq 1,\ i \in I_{n}.
	$$
  There exists a coupling of $(\tau_{i}^{n-1})_{1 \le i \le n-1}$ and $(\tau_{i}^{n})_{1 \le i \le n}$ such that 
  	$$
	\tau_{i-1}^{n-1} \leq \tau_{i}^{n},\ 2 \leq i \leq n.
	$$
  \end{lemma}
  \begin{proof}
    Let $Y_{2},\dots, Y_{n}$ be iid random variables with exponential
    distribution $\mathcal{E}(1)$. Set $\tau_{1}^{n} \coloneq \tau_{1}^{n-1} \coloneq 1$ and
	$$
	\tau_{i}^{n} \coloneq \Big\lceil -\frac{Y_{i}}{\log(1-p_{i}^{n})} \Big\rceil\ \text{and\ } \tau_{i-1}^{n-1} \coloneq \Big\lceil -\frac{Y_{i}}{\log(1-p_{i-1}^{n-1})} \Big\rceil,\ 2 \leq i \leq n
	$$
	Then $\tau_{i}^{n} \sim \operatorname{Geom}(p_{i}^{n})$ and
	$\tau_{i}^{n-1} \sim \operatorname{Geom}(p_{i}^{n-1})$. This is simply a consequence of the fact that if $Y \sim \mathcal{E}(1)$,
	then:
	$$
	\Big\lceil \frac{Y}{\lambda} \Big\rceil \sim \operatorname{Geom}(e^{-\lambda}),\ \lambda > 0.
	$$
	The inequality of the lemma is proved by noticing that
	$p\mapsto \lceil -y/\log(1-p) \rceil$ is non-increasing on $(0,1)$ for all
	$y\in \R_{+}$ and that $p_{i}^{n} \leq p_{i-1}^{n-1}$:
	$$
	p_{i}^{n} = \frac{n-i+1}{n} \leq \frac{n-i+1}{n-1} = p_{i-1}^{n-1}.
	$$
	This concludes the proof.
  \end{proof}
  As a consequence,  we can compare $Z_{n}$ and $Z_{n-1}$
  in $\LL^{1}$-norm:
  \begin{lemma}
\label{lem_coupon_core:2}
	There exists a positive constant $C$ such that:
	$$
	\Vert Z_{n} - Z_{n-1} \Vert_{\LL^{1}} \leq C\ \frac{\log n}{n},\ n\geq 2.
	$$
  \end{lemma}

  \begin{proof}
	Notice that since $\tau_{1}^{n} = 1$ a.s., we have for all $n \geq 2$:
	\begin{align*}
	  Z_{n} - Z_{n-1} & = \frac{1}{n}\sum_{i=1}^{n}\tau_{i}^{n} - \log n - \frac{1}{n-1}\sum_{i=1}^{n-1}\tau_{i}^{n-1} + \log(n-1)                                                   \\
					  & = \frac{1}{n}\sum_{i=1}^{n-1}\tau_{i+1}^{n} - \frac{1}{n-1}\sum_{i=1}^{n-1}\tau_{i}^{n-1} + \frac{1}{n} + \log\Big(1 - \frac{1}{n}\Big)                      \\
					  & = \frac{1}{n}\sum_{i=1}^{n-1}(\tau_{i+1}^{n} - \tau_{i}^{n-1})- \frac{1}{n(n-1)}\sum_{i=1}^{n-1}\tau_{i}^{n-1} + \frac{1}{n} + \log\Big(1 - \frac{1}{n}\Big).
	\end{align*}
	Since $\tau_{i+1}^{n} \geq \tau_{i}^{n-1}$ a.s. and $\vert \log(1 - n^{-1}) \vert \leq (n-1)^{-1}$, we have:
	\begin{align*}
	  \Vert Z_{n} - Z_{n-1} \Vert_{\LL^{1}} & \leq \frac{2}{n-1} + \frac{1}{n}\sum_{i=1}^{n-1}\Vert \tau_{i+1}^{n} - \tau_{i}^{n-1} \Vert_{\LL^{1}} + \frac{1}{n(n-1)}\sum_{i=1}^{n-1}\Vert \tau_{i}^{n-1} \Vert_{\LL^{1}} \\
												   & = \frac{2}{n-1} + \frac{1}{n}\sum_{i=1}^{n-1}\esp[\tau_{i+1}^{n} - \tau_{i}^{n-1}] + \frac{1}{n(n-1)}\sum_{i=1}^{n-1}\esp[\tau_{i}^{n-1}]                                                  \\
												   & = \frac{2}{n-1} + \frac{2}{n}H_{n-1},
	\end{align*}
	where $H_{n} \coloneq \sum_{i=1}^{n}i^{-1}$ denotes the $n$-th harmonic
	number.
  \end{proof}
  \begin{lemma}
	\label{lem_coupon_core:1}
	The density of the distribution of  $(\tau_{1}^{n},\dots,\tau_{n-1}^{n})$ with
	respect to the law of $(\tau_{1}^{n-1},\dots,\tau_{n-1}^{n-1})$ is
	equal to:
	\begin{align*}
	  \frac{P_{n-1}^{n}(t_{1},\dots,t_{n-1})}{P_{n-1}^{n-1}(t_{1},\dots,t_{n-1})} & =  K_{n}\,\Big(1-\frac{1}{n}\Big)^{(n-1)\big(\frac{1}{n-1}\sum_{i=1}^{n-1}(t_{i}-\log(n-1))\big)},
	\end{align*}
	where $P_{j}^{n}$ denotes the joint
	distribution of $(\tau_{1}^{n},\dots,\tau_{j}^{n})$ for all $j \in I_{n}$.
  \end{lemma}
  \begin{proof}
	For any $(t_{1},\cdots,t_{n})\in \N^{*}$,  we have
	\begin{multline*}
	  \frac{P_{n-1}^{n}(t_{1},\dots,t_{n-1})}{P_{n-1}^{n-1}(t_{1},\dots,t_{n-1})} \\
	  \begin{aligned}
	& = \frac{\prod_{i=1}^{n-1}\big(\frac{n-i+1}{n}\big)\big(\frac{i-1}{n}\big)^{t_{i}-1}}{\prod_{i=1}^{n-1}\big(\frac{n-i}{n-1}\big)\big(\frac{i-1}{n-1}\big)^{t_{i}-1}} \\
	& = \Big(1-\frac{1}{n}\Big)^{n-1}\frac{n!}{(n-1)!}\frac{\prod_{i=1}^{n-1}\big(\frac{1}{n}\big)^{t_{i}-1}}{\prod_{i=1}^{n-1}\big(\frac{1}{n-1}\big)^{t_{i}-1}}         \\
		  & = n\Big(1-\frac{1}{n}\Big)^{\sum_{i=1}^{n-1}t_{i}}                                                                                                                  \\
  & = n\Big(1-\frac{1}{n}\Big)^{(n-1)\big(\frac{1}{n-1}\sum_{i=1}^{n-1}t_{i}\big)}                                                                                      \\
	  & = n\Big(1-\frac{1}{n}\Big)^{(n-1)\log(n-1)}\Big(1-\frac{1}{n}\Big)^{(n-1)\big(\frac{1}{n-1}\sum_{i=1}^{n-1}(t_{i}-\log(n-1))\big)},
	  \end{aligned}
	\end{multline*}
	hence the announced result.
  \end{proof}
  \begin{proof}[Proof of Theorem~\ref{thm_coupon_core:main_identity}]

	The idea is to do a change of probability: we replace the distribution of
	$(\tau_{1}^{n},\dots,\tau_{n-1}^{n},\tau_{n}^{n})$ by the distribution of
	$(\tau_{1}^{n-1},\dots,\tau_{n-1}^{n-1},\tau_{n}^{n})$. The density function of the former 	with respect to the latter will behave as $e^{-Z_{n}} \simeq e^{-Z}$
	when $n$ goes to infinity, where $Z \sim \mathcal{G}(0,1)$. As for
	$\tau_{n}^{n}$, this term is singled out because of the following well-known
	result:
	$$
	G_{n}\coloneq\frac{\tau_{n}^{n}}{n} \eqdis \frac{1}{n}\operatorname{Geom}\Big(\frac{1}{n}\Big) \xrightarrow[n\to\infty]{\text{law}} \mathcal{E}(1).
	$$
	This will give us the exponential random variable that appears in the
	generator of the Gumbel distribution. 
	\begin{multline*}
\esp\big[f'(Z_{n})\big]  = \esp\Big[f'\Big(\frac{1}{n}\sum_{i=1}^{n-1}\tau_{i}^{n}-\log n + \frac{\tau_{n}^{n}}{n}\Big)\Big]                                                                                                            \\
\begin{aligned}
& = \esp\Big[f'\Big(\frac{1}{n}\sum_{i=1}^{n-1}\tau_{i}^{n}-\log n + G_{n}\Big)\Big]                                                                                                                   \\
							  & = \esp\Big[\frac{P_{n-1}^{n}(\tau_{1}^{n-1},\dots,\tau_{n-1}^{n-1})}{P_{n-1}^{n-1}(\tau_{1}^{n-1},\dots,\tau_{n-1}^{n-1})} \ f'\big(\frac{1}{n}\sum_{i=1}^{n-1}\tau_{i}^{n-1}-\log n + G_{n}\big)\Big] \\
							  & = K_{n}\ \esp\Big[\Big(1-\frac{1}{n}\Big)^{(n-1)Z_{n-1}}f'\big(\frac{1}{n}\sum_{i=1}^{n-1}\tau_{i}^{n-1}-\log n + G_{n}\big)\Big]                                                                      \\
& = K_{n}\,\esp\big[e^{-\beta_{n}Z_{n-1}}f'\big(\alpha_{n}Z_{n-1} + {G_{n}} + \delta_{n}\big)\big].
\end{aligned}
\end{multline*}
The proof is thus complete.
  \end{proof}
  \begin{notation}
    The following notation will be convenient in the sequel:
  \begin{align}\label{cpn_g_n}
	g_{n}(x) \coloneq e^{-\beta_{n}x}\esp\big[f'\big(\alpha_{n}x + G_{n}+ \delta_{n}\big)\big].
  \end{align}
  \end{notation}
  We can then proceed to the main part of the proof.
  \begin{theorem}
\label{cpn_aroundI} Let $h \in \lipbis$. There exists a
	constant $C>0$  such that:
	$$
	\big\vert\esp\big[\Lo \Po_{t}h(Z_{n})\big]\big\vert \leq C\ \frac{1 + \vert\log \gamma_{t}\vert}{\gamma_{t}}\frac{\log n}{n},\ \forall t \in [1, \infty),\ \forall n \ge 1.
	$$
\end{theorem}
On the other hand, we have the following bound for $\esp[\Lo \Po_{t}h(Z_{n})]$
in the neighborhood of~ $0$:

\begin{theorem}
\label{cpn_around0} Let $h\in \lipbis$. There exists a constant $C>0$ such that:
	$$
	\vert\esp[\Lo\Po_{t}h(Z_{n})]\vert \leq C \ \frac{\log n}{n},\ \forall t \in [0, 1 ],\ \forall n \ge 1.
	$$
\end{theorem}
The proofs of these two theorems are postponed to the last section. Taking them
for granted,  we can then conclude and obtain Theorem~\Ref{cpn_main}.
 \begin{proof}[Proof of Theorem~\protect\ref{cpn_main}]
	We already know that for every $h \in \lipbis$
	$$
	\esp[h(Z_{n})] - \esp[h(Z)] = -\int_{0}^{\infty}\esp\big[\Lo \Po_{t}h(Z_{n})\big]\dif t.
  $$
In view of Theorems \ref{cpn_aroundI}
	and \ref{cpn_around0}, we have:
	\begin{align*}
	  \vert \esp[h(Z_{n})] - \esp[h(Z)] \vert & \leq \int_{0}^{\infty}\big\vert\esp\big[\Lo \Po_{t}h(Z_{n})\big]\big\vert\dif t                                                                                                                           \\
											  & = \int_{0}^{1} \big\vert\esp\big[\Lo \Po_{t}h(Z_{n})\big]\big\vert\dif t + \int_{1}^{\infty} \big\vert\esp\big[\Lo \Po_{t}h(Z_{n})\big]\big\vert\dif t \\
											  & \le C\ \frac{\log n}{n}\Big(1 + \int_{1}^{\infty}\frac{1 + \vert \log \gamma_{t}  \vert}{\gamma_{t}}\dif t\Big)\cdotp
	\end{align*}
  It remains to remark that the rightmost integral is finite.
  \end{proof}

\section{Proofs of theorems \ref{cpn_aroundI} and \ref{cpn_around0}}

\subsection{Proof of theorem \ref{cpn_aroundI} - Bounding the error on $[1, +\infty)$}
We introduce the decomposition:
\begin{multline}
\label{cpn_dcptn}
  \esp[  \Lo f(Z_{n})]\\
	\shoveleft{ = -\esp[h'(Z_{n})] + \esp\big[e^{-Z_{n}}h'(Z_{n} + Y)\big]}\\
	\shoveleft{	 = -K_{n}\esp[g_{n}(Z_{n-1})] + \esp[g(Z_{n})]}                                                                                     \\
	\shoveleft{= (1-K_{n})\esp[g_{n}(Z_{n-1})] + \big(\esp[g(Z_{n})] - \esp[g(Z_{n-1})]\big)} \\
	\shoveright{+ \big(\esp[g(Z_{n-1})]
	- \esp[g_{n}(Z_{n-1})]\big)}\\
		 = A_1 + A_2 + A_3,
\end{multline}
where the second equality is a consequence of our Stein characterization of the Gumbel distribution \eqref{cpn_Stein_lemma}. The function $g_{n}$ has been defined in \eqref{cpn_g_n}, with $f = \Po_{t}h$, so that, thanks to identity \eqref{P_tDiff}:
$$
g_{n}(x) = e^{-\beta_{n}x}\esp[e^{-\gamma_{t}e^{-(\alpha_{n}x + G_{n} + \delta_{n})}}h'(\alpha_{n}x + G_{n} + \delta_{n} - t)],
$$
and $g$ is the pointwise limit of $g_{n}$:
\begin{align}\label{def_g}
	g(x) \coloneq e^{-x}\esp[e^{-\gamma_{t}e^{-(x + Y)}}h'(x + Y - t)],
\end{align}
with $Y \sim \mathcal{E}(1)$.

We deal first with $A_1$ and $A_2$.
\begin{lemma}\label{lemma_A12} Let $h\in \lipbis$. We couple $G_{n}$ with $Y$ in the following
	way
  \begin{equation}
    \label{eq_coupon_core:4}
    	G_{n} \coloneq \frac{1}{n}\Big\lceil \frac{-Y}{\log \alpha_{n}} \Big\rceil \sim \frac{1}{n}\operatorname{Geom}\Big(\frac{1}{n}\Big).
  \end{equation}
	Then, the function $g$ is $(e\gamma_{t})^{-1}$-Lipschitz and
	 for $t > 0$, $n \ge 1$, $i=1,2$, we have
\begin{equation}
	\label{eq:coupon_core:5}
					A_i \leq
					\frac{C}{\gamma_{t}}\frac{\log n}{n}\cdotp
\end{equation}
\end{lemma}
\begin{proof}
Note that $g'$ admits the following simple expression
 thanks to the change of variable $z = y+x$:
 \begin{align*}
   g(x) &= e^{-x}\int_{0}^{\infty}e^{-\gamma_{t}e^{-(x+y)}}h'(x+y-t)e^{-y}\dif y \\
   &= \int_{x}^{\infty}e^{-z}e^{-\gamma_{t}e^{-z}}h'(z-t)\dif z.
\end{align*}
			Differentiating with respect to $x$, we get $g'(x) = -e^{-x}e^{-\gamma_{t}e^{-x}}h'(x-t)$, so that
		      $$
			\vert g'(x) \vert \leq e^{-x}e^{-\gamma_{t}e^{-x}} \leq \frac{e^{-1}}{\gamma_{t}}\cdotp
		      $$
Combined with Lemma~\ref{lem_coupon_core:2}, this yields~\eqref{eq:coupon_core:5} for $A_{2}$.
	 Recall that $g_{n}$ is given by \eqref{cpn_g_n}, and that the
			constants $\alpha_{n}, \beta_{n}$ and $\delta_{n}$ have been defined
			just before \eqref{cpn_ineq_cstt}. We know that the double
			exponential term is a non-decreasing function, and that
			$\lceil y \rceil \leq y+1$, for every non-negative $y$. Besides, as
			$-n\log\alpha_{n} \geq 1$ and $\delta_{n} \leq -1/n$, we have:
			\begin{align}\label{cpn_Y_G_n}
			  - \frac{Y}{n\log\alpha_{n}} + \delta_{n} \leq {G_{n}} + \delta_{n} \leq - \frac{Y}{n\log\alpha_{n}} + \delta_{n} + \frac{1}{n} \leq Y.
			\end{align}
			Furthermore, for $x\ge -\log n$, we have
			\begin{align*}
			  \vert g_{n}(x) \vert \leq e^{-\beta_{n}x}\esp\Big[e^{-\gamma_{t}e^{-(\alpha_{n}x + G_n  + \delta_{n})}}\Big] & \leq e^{-\beta_{n}x}\esp\Big[e^{-\gamma_{t}e^{-(\alpha_{n}x + Y)}}\Big]                            \\
								   & = \frac{1}{\gamma_{t}}e^{-(\beta_{n}-\alpha_{n})x}\big(1-e^{-\gamma_{t}e^{-\alpha_{n}x}}\big)      \\
								   & \leq \frac{1}{\gamma_{t}}e^{-(\beta_{n}-\alpha_{n})x} \leq \frac{1}{\gamma_{t}}n^{\beta_{n}-\alpha_{n}} \leq \frac{C}{\gamma_{t}},
			\end{align*}
			since
			$n^{\beta_{n}-\alpha_{n}}$ tends to $1$ as $n$ goes to infinity.
Moreover, it holds that
			\begin{align}\label{K_n}
0 \leq K_{n}-1&\leq ne^{-(1-\frac{1}{n})\log(n-1)}-1 
			 \\ & \underset{n\to\infty}{\sim} \frac{\log n}{n},
\end{align}
where $a_{n} \underset{n \to \infty}{\sim} b_{n}$ means that $a_{n}/b_{n}$ converges to $1$ when $n$ goes to infinity. 
Finally, remark that $Z_{n-1} \geq 1-\log(n-1) \geq -\log
n$ a.s. to conclude that~\eqref{eq:coupon_core:5} holds for $A_1$.
\end{proof}
\begin{remark}
The intermediate result, available for $x$ greater than $-\log n$
			\begin{equation}\label{cpn_bound_g_n}
			  e^{-\beta_{n}x}\esp\Big[e^{-\gamma_{t}e^{-(\alpha_{n}x + G_n  + \delta_{n})}}\Big] \leq \frac{C}{\gamma_{t}},
			\end{equation}
			will be useful on its own in the sequel.
    \end{remark}
As for $A_3$, we will need to decompose it further:
\begin{multline*}
  \esp[g(Z_{n-1})] - \esp\big[g_{n}(Z_{n-1})\big] \\
  \shoveleft{	= \Big(\esp[g(Z_{n-1})] - \esp\big[e^{-\beta_{n}Z_{n-1}}e^{-\gamma_{t}e^{-(\alpha_{n}Z_{n-1} + G_n  + \delta_{n})}}h'(Z_{n-1} + Y - t)\big]\Big)}                   \\
  + \Big(\esp\big[e^{-\beta_{n}Z_{n-1}}e^{-\gamma_{t}e^{-(\alpha_{n}Z_{n-1} + G_n  + \delta_{n})}}h'(Z_{n-1} + Y - t)\big] - \esp\big[g_{n}(Z_{n-1})\big]\Big) \\
  =B_{1}+B_{2}.
\end{multline*}
\begin{lemma}
  \label{lem_coupon_core:3}
  Let $h \in \lipbis$. There exists $C>0$ such that for any $n\ge 1$, for any $t > 0$,
  \begin{equation*}
|B_{2}|\le \frac{C}{n\gamma_{t}}\cdotp
  \end{equation*}
\end{lemma}
We begin with two technical lemmas.
\begin{lemma}\label{cpn_exp_bnd}
  For all $\lambda \geq 0$, the sequence $(\esp[e^{-\lambda Z_{n}}])_{n\geq 1}$
  is bounded. Moreover, we have the asymptotic expansion:
  \begin{align}\label{cpn_Z_n_asym}
	\esp[e^{-Z_{n}}] - 1 \underset{n \to \infty}{\sim} - \frac{\log n}{2n}.
  \end{align}
\end{lemma}

\begin{proof} For all $n\geq 1$, we have by definition of $Z_{n}$:
  \begin{align*}
    \esp[e^{-\lambda Z_{n}}] = n^{\lambda}\prod_{i=1}^{n}\frac{1 - \frac{i-1}{n}}{e^{\frac{\lambda}{n}} - \frac{i-1}{n}} &= n^{\lambda}n!\prod_{i=1}^{n}\frac{1}{ne^{\frac{\lambda}{n}} - i + 1}\\
                             &\leq n^{\lambda}n!\prod_{i=1}^{n}\frac{1}{n - i + 1 + \lambda} = n^{\lambda}n!\prod_{i=1}^{n}\frac{1}{i+\lambda},
\end{align*}
  thanks to the inequality $ne^{\lambda/n} \geq n+\lambda$. Therefore, thanks to the identity $x = \exp(\log x)$ and the inequality $\log i - \log(i + \lambda) \leq -\lambda/(i + \lambda)$, we have:
  \begin{align*}
  \sup_{n\ge 1}  \esp[e^{-\lambda Z_{n}}] & \leq \sup_{n\ge 1}e^{\lambda \log n + \sum_{i=1}^{n}\log i - \log(i+\lambda)}\\
  										 & \leq \sup_{n\ge 1}e^{-\lambda\big(\sum_{i=1}^{n}\frac{1}{i+\lambda} - \log n\big)}\\
                                          &= \sup_{n\ge 1} e^{-\lambda(H_{n} - \log n)} e^{-\lambda\sum_{i=1}^{\infty}\frac{1}{i(i+\lambda)}}<\infty.
\end{align*}
  As for the second part of the
  statement
  \begin{align*}
    \esp[e^{-Z_{n}}] = nn!\prod_{i=1}^{n}\frac{1}{ne^{1/n} - i + 1} & = \frac{nn!}{(n+1)!}\prod_{i=1}^{n}\frac{n - i + 2}{ne^{1/n} - i + 1}                              \\
					 & = \frac{n}{n+1}\prod_{i=1}^{n}\frac{n - i + 2}{ne^{1/n} - i + 1}                                   \\
					 & = \frac{n}{n+1}\prod_{i=1}^{n}\frac{n - i + 2}{(n(e^{1/n}-1) - 1) + n - i + 2}                     \\
					 & = \frac{n}{n+1}\prod_{i=1}^{n}\Big(1 + \frac{n(e^{1/n} - 1) - 1}{i + 1}\Big)^{-1}                  \\
					 & = \frac{n}{n+1}\exp\Big(-\sum_{i=1}^{n}\log\big(1 + \frac{n(e^{1/n} - 1) - 1}{n - i + 2}\big)\Big).
  \end{align*}
	Therefore, a standard result about positive divergent series yields:
  \begin{align*}
    \esp[e^{-Z_{n}}] - 1 &\underset{n\to \infty}{\sim} \frac{n}{n+1}e^{-\frac{1}{2n}H_{n}} - 1\\
 & \underset{n\to \infty}{\sim} -\frac{1}{2n}\log n,
  \end{align*}
	which concludes the proof.
  \end{proof}
\begin{lemma}
  \label{lem_coupon_core:4}
  Let $Y \sim \mathcal{E}(1)$ and $G_{n}$ defined as in~\eqref{eq_coupon_core:4}. Then,
  we have :
  \begin{equation*}
    \Vert Y - G_n  \Vert_{\LL^{2}}\le \frac{1}{n} , n\ge 2.
  \end{equation*}
  \end{lemma}
\begin{proof}
  Set
  $u_{n} \coloneq \int_{0}^{n\lambda}ye^{-y}\dif y - 1$.
  We have
  \begin{align*}
    \esp\Big[Y\big\lceil \frac{Y}{\lambda} \big\rceil\Big] = \int_{0}^{\infty}y\big\lceil \frac{y}{\lambda} \big\rceil e^{-y}\dif y & = \sum_{n=0}^{\infty}(n+1)\int_{n\lambda}^{(n+1)\lambda}ye^{-y}\dif y                                 \\
														   & = \sum_{n=0}^{\infty}(n+1)(u_{n+1} - u_{n})                                                           \\
														   & = -\sum_{n=0}^{\infty}u_{n} = \frac{\lambda e^{-\lambda}}{(1-e^{-\lambda})^{2}} + \frac{1}{1-e^{-\lambda}}\cdotp
  \end{align*}
  Evaluating this identity at $\lambda = -\log(1-1/n) > 0$, we get that
  \begin{align*}
	\esp\big[\vert Y - G_n  \vert^{2}\big] & = 4 - \frac{1}{n} -\frac{2}{n}\esp\Big[Y\big\lceil \frac{Y}{\lambda} \big\rceil\Big] \\
															  & = 2 - \frac{1}{n} + 2(n-1)\log\Big(1-\frac{1}{n}\Big)                                    \\
															  & \leq \frac{1}{n^{2}}\cdotp
  \end{align*}
  The last line stems from the well-known inequality:
	$$ \log\Big(1-\frac{1}{x}\Big) \leq -\frac{1}{x} - \frac{1}{2x^{2}},\ x > 1.
	$$
  The proof is thus complete.
  \end{proof}

  \begin{proof}[Proof of Lemma~\protect~\ref{lem_coupon_core:3}]
    For the sake of notations, we set temporarily
    \begin{equation*}
\zeta_{n}=e^{-\beta_{n}Z_{n-1}}e^{-\gamma_{t}e^{-(\alpha_{n}Z_{n-1} + G_n  + \delta_{n})}}.
    \end{equation*}
    Recall that $\delta_{n}$ is non-positive, hence:
\begin{multline*}
  \Big\vert\esp\Big[e^{-\beta_{n}Z_{n-1}}e^{-\gamma_{t}e^{-(\alpha_{n}Z_{n-1} + G_n  + \delta_{n})}} h'(Z_{n-1} + Y - t)\Big] - \esp\big[g_{n}(Z_{n-1})\big]\Big\vert                                                                                                                           \\
  \begin{aligned}
  &\leq \esp\Big[\zeta_{n}\big\vert h'(Z_{n-1} + Y - t) - h'(\alpha_{n}Z_{n-1} + G_n  + \delta_{n} - t) \big\vert \Big]                                         \\
  & \leq \esp\Big[\zeta_{n}\big\vert (1-\alpha_{n})Z_{n-1} + Y - G_n  - \delta_{n}\big\vert \Big]                                                                 \\
 & \leq \esp\Big[\zeta_{n}\Big(\frac{1}{n}\vert Z_{n-1} \vert + \vert Y - G_n \vert - \delta_{n}\Big)\Big].
  \end{aligned}
\end{multline*}
That last expression is also equal to:
\begin{multline*}
  \frac{1}{n}\esp\Big[e^{-\beta_{n}Z_{n-1}}e^{-\gamma_{t}e^{-(\alpha_{n}Z_{n-1} + G_n  + \delta_{n})}}(\vert Z_{n-1} \vert - n\delta_{n})\Big]\\
  \shoveright{ + \esp\Big[e^{-\beta_{n}Z_{n-1}}e^{-\gamma_{t}e^{-(\alpha_{n}Z_{n-1} + G_n  - \delta_{n})}}\vert Y - G_n \vert\Big]   }                                                                                                        \\
  \begin{aligned}
    & \leq \frac{C}{n\gamma_{t}}\big(\esp[\vert Z_{n-1} \vert] - n\delta_{n}\big) + \Big(\esp\big[e^{-2\beta_{n}Z_{n-1}}e^{-2\gamma_{t}e^{-(\alpha_{n}x + G_n  + \delta_{n})}}\big]\Big)^{\frac{1}{2}}\Vert Y - G_n \Vert_{\LL^{2}}\\
    & \leq \frac{C}{n\gamma_{t}}\big(\esp[\vert Z_{n-1} \vert] - n\delta_{n}\big) + \frac{C}{\gamma_{t}}\big(\esp[e^{-2\beta_{n}Z_{n-1}}]\big)^{\frac{1}{2}}\Vert Y - G_n \Vert_{\LL^{2}},
  \end{aligned}
\end{multline*}
according to inequality~\eqref{cpn_bound_g_n} and Cauchy-Schwarz
inequality.  It is well known that
\begin{equation*}
  \sup_{n\ge 1}\esp[Z_{n}^{2}]<\infty \text{ hence } \sup_{n\ge 1}\esp[|Z_{n}|]<\infty.
\end{equation*}
The result then follows from Lemma~\ref{cpn_exp_bnd} and Lemma~\ref{lem_coupon_core:4}.
\end{proof}
 Recall that
\begin{multline*}
  B_{1}=\esp[e^{-Z_{n-1}}e^{-\gamma_{t}e^{-(Z_{n-1} + Y)}}h'(Z_{n-1} + Y - t)]\\ - \esp\big[e^{-\beta_{n}Z_{n-1}}e^{-\gamma_{t}e^{-(\alpha_{n}Z_{n-1} + G_n  + \delta_{n})}}h'(Z_{n-1} + Y - t)\big] = C_{1} + C_{2},
\end{multline*}
where
\begin{align*}
  C_{1}&=\esp[e^{-Z_{n-1}}(e^{-\gamma_{t}e^{-(Z_{n-1} + Y)}}-e^{-\gamma_{t}e^{-(\alpha_{n}Z_{n-1} + G_n  + \delta_{n})}})h'(Z_{n-1} + Y - t)],\\
  C_{2}&= \esp\Big[\big(e^{-Z_{n-1}} - e^{-\beta_{n}Z_{n-1}}\big)e^{-\gamma_{t}e^{-(\alpha_{n}Z_{n-1} + G_n  + \delta_{n})}}  h'(Z_{n-1} + Y - t)\Big].
\end{align*}
\begin{lemma}\label{lemma_C2}
  There exists $C>0$ such that
  \begin{equation}
    \label{eq_coupon_core:11}
    |C_{2}|\le\frac{C}{n\gamma_{t}} \text{ for } n\ge 2.
  \end{equation}
\end{lemma}
\begin{proof}
  We have for any $h \in \lipbis$
  \begin{align*}
    |C_{2}|  & \leq  \esp\Big[\big\vert e^{-\beta_{n}Z_{n-1}} - e^{-Z_{n-1}}\big\vert e^{-\gamma_{t}e^{-(\alpha_{n}Z_{n-1} + G_n  + \delta_{n})}}\Big] \\
	  & \leq \esp\Big[\big\vert e^{-\beta_{n}Z_{n-1}} - e^{-Z_{n-1}}\big\vert e^{-\gamma_{t}e^{-(\alpha_{n}Z_{n-1} + Y)}}\Big]                             \\
	  & = \esp\Big[\vert e^{-\beta_{n}Z_{n-1}} - e^{-Z_{n-1}}\vert\varphi_{n}(Z_{n-1})\Big],
  \end{align*}
  where we set
  $\varphi_{n}(x) \coloneq \esp\big[e^{-\gamma_{t}e^{-(\alpha_{n}x + Y)}}\big]$
  for all $x \geq -\log n$. The change of variable $z = y - \log \gamma_{t}$ shows that:
$$
\varphi_{n}(x) = \frac{1}{\gamma_{t}}\int_{-\log \gamma_{t}}^{\infty}e^{-z}e^{-e^{-(\alpha_{n}x + z )}}\dif z \leq \frac{1}{\gamma_{t}}\int_{\R}e^{-z}e^{-e^{-(\alpha_{n}x + z)}}\dif z = \frac{1}{\gamma_{t}}e^{\alpha_{n}x}.
$$
According to Cauchy-Schwarz inequality and Lemma~\ref{cpn_exp_bnd}, it follows
that
\begin{multline*}
  \Big\vert\esp\Big[(e^{-\beta_{n}Z_{n-1}} - e^{-Z_{n-1}})\varphi_{n}(Z_{n-1})\big]\big\vert \\
  \begin{aligned}
  &\leq \frac{C}{\gamma_{t}}\esp\big[\big\vert e^{-(\beta_{n}-\alpha_{n})Z_{n-1}} - e^{-(1-\alpha_{n})Z_{n-1}}\Big\vert\Big]                            \\
  &\leq \frac{C}{\gamma_{t}}(1-\beta_{n})\esp\Big[\vert Z_{n-1} \vert \big(e^{-(\beta_{n} - \alpha_{n})Z_{n-1}} + e^{-(1-\alpha_{n})Z_{n-1}}\big)\Big]\\
  &\le \frac{C}{n\gamma_{t}}\cdotp
  \end{aligned}
\end{multline*}
The proof is this complete.
\end{proof}
Before bounding $C_{1}$, we need another technical lemma.
\begin{lemma}
  \label{lem_coupon_core:6}
  There exists a constant $C > 0$ such that for all $n\geq 2$, we have for every $t > 0$ and $x \in \R$:
  \begin{multline*}
	e^{-x}\Big\vert \int_{0}^{\infty}e^{-y}e^{-\gamma_{t}e^{-(\alpha_{n}x + y)}}\dif y - \varepsilon_{n}\int_{0}^{\infty}e^{-\varepsilon_{n} y}e^{-\gamma_{t}e^{-(\alpha_{n}x + y + \delta_{n})}} \dif y \Big\vert                                                                                     \\		                                                                                                                                                                     \leq C\, \frac{\vert x \vert + \vert \log\gamma_{t}\vert}{\gamma_{t}}\frac{\log n}{n}e^{-x/n}.
  \end{multline*}
  In particular
  \begin{multline*}
	\esp\Big[e^{-Z_{n-1}}\Big\vert \int_{0}^{\infty}e^{-y}e^{-\gamma_{t}e^{-(\alpha_{n}Z_{n-1} + y)}}\dif y - \varepsilon_{n}\int_{0}^{\infty}e^{-\varepsilon_{n} y}e^{-\gamma_{t}e^{-(\alpha_{n}Z_{n-1} + y + \delta_{n})}} \dif y \Big\vert\Big]                                                    \\
	\leq C\frac{1 + \vert \log \gamma_{t}\vert}{\gamma_{t}}\frac{\log n}{n}\cdotp
  \end{multline*}
\end{lemma}
\begin{proof}
  The arguments are similar to those given above, but certain new difficulties
  arise nonetheless.

 First we replace $\varepsilon_{n} = -n\log(1-1/n) \geq 1$ by 1. Recall that $\delta_{n}$ is negative, so
		  that:
		  \begin{multline*}
			0 \leq (\varepsilon_{n}-1)\int_{0}^{\infty}e^{-\varepsilon_{n} y}e^{-\gamma_{t}e^{-(\alpha_{n}x + y + \delta_{n})}}\dif y\\ \leq (\varepsilon_{n}-1)\int_{0}^{\infty}e^{-y}e^{-\gamma_{t}e^{-(\alpha_{n}x + y)}}\dif y \leq \frac{e^{\alpha_{n}x}}{\gamma_{t}}(\varepsilon_{n}-1).
		  \end{multline*}
		  Since, $z \mapsto e^{-yz}$ is Lipschitz on
		  $[1,\varepsilon_{n}]$, we have:
		  \begin{align*}
			0 \leq \int_{0}^{\infty}(e^{-y} - e^{-\varepsilon_{n} y}) & e^{-\gamma_{t}e^{-(\alpha_{n}x + y + \delta_{n})}}\dif y                                                                                                                                       \\
																	  & \leq (\varepsilon_{n}-1)\int_{0}^{\infty}ye^{-y}e^{-\gamma_{t}e^{-(\alpha_{n}x + y + \delta_{n})}}\dif y\\
																	  &= (\varepsilon_{n}-1)e^{-\tau}\esp\big[(Z + \tau)\ind_{\lb Z + \tau \geq 0 \rb}\big].
		  \end{align*}
		  where $\tau \coloneq \log \gamma_{t} -\alpha_{n}x - \delta_{n} \in \R$
		  and $Z \sim \mathcal{G}(0,1)$.
Thus:
		  \begin{multline*}
			0 \leq \int_{0}^{\infty}(e^{-y} - e^{-\varepsilon_{n} y}  )e^{-\gamma_{t}e^{-(\alpha_{n}x + y + \delta_{n})}}\dif y  \\
			\begin{aligned}
			  & \leq e^{-\tau}\big(C + \vert\log \gamma_{t}\vert -\alpha_{n}\vert x \vert\big)(\varepsilon_{n}-1)\\
			  &= \frac{e^{\alpha_{n}x + \delta_{n}}}{\gamma_{t}}\big(C + \vert\log \gamma_{t}\vert + \alpha_{n}\vert x \vert\big)(\varepsilon_{n}-1).
			\end{aligned}
		  \end{multline*}
 As for the final term, because
		  $\delta_{n} = \log(1-1/n) - \log(n-1)/n \leq 0$, we find:
		      \begin{multline*}
				0 \leq \int_{0}^{\infty}e^{-\varepsilon_{n}y}e^{-\gamma_{t}e^{-(\alpha_{n}x + y)}}\dif y -  \int_{0}^{\infty}e^{-\varepsilon_{n} y}e^{-\gamma_{t}e^{-(\alpha_{n}x + y + \delta_{n})}}\dif y                                                                                                                \\
				\begin{aligned}
				  & = \int_{0}^{\infty}e^{-\varepsilon_{n}y}e^{-\gamma_{t}e^{-(\alpha_{n}x + y)}}\dif y - e^{\varepsilon_{n}\delta_{n}}\int_{\delta_{n}}^{\infty}e^{-\varepsilon_{n} y}e^{-\gamma_{t}e^{-(\alpha_{n}x + y)}}\dif y \\
				  & \leq \int_{0}^{\infty}e^{-\varepsilon_{n}y}e^{-\gamma_{t}e^{-(\alpha_{n}x + y)}}\dif y - e^{\varepsilon_{n}\delta_{n}}\int_{0}^{\infty}e^{-\varepsilon_{n} y}e^{-\gamma_{t}e^{-(\alpha_{n}x + y)}}\dif y       \\
				  & = (1-e^{\varepsilon_{n}\delta_{n}})\int_{0}^{\infty}e^{-\varepsilon_{n}y}e^{-\gamma_{t}e^{-(\alpha_{n}x + y)}}\dif y                                                                                           \\
				  & \leq \frac{e^{\alpha_{n}x}}{\gamma_{t}}(1-e^{\varepsilon_{n}\delta_{n}}).
				\end{aligned}
			  \end{multline*}
	The second inequality of the statement is obtained by replacing $x$ by
	$Z_{n-1}$ and integrating with respect to this random variable.
  \end{proof}
  \begin{lemma}
    \label{lem_coupon_core:5}
    With the notations of Theorem~\ref{cpn_aroundI}, there exists $C>0$ such
    that
    \begin{equation*}
      |C_{1}|\le C \frac{1+\log \gamma_{t}}{\gamma_{t}} \frac{\log n}{n}\text{
        for } t\in [1,\infty), \ n\ge 1.
    \end{equation*}
  \end{lemma}
\begin{proof}
We have
  \begin{multline}
\label{eq_coupon_core:12}|C_{1}|\le   \esp\Big[e^{-Z_{n-1}}\big\vert e^{-\gamma_{t}e^{-(Z_{n-1} + Y)}} - e^{-\gamma_{t}e^{-(\alpha_{n}Z_{n-1} + Y)}}\big\vert\Big] \\
  \shoveright{	+ \esp\Big[e^{-Z_{n-1}}\big( e^{-\gamma_{t}e^{-(\alpha_{n}Z_{n-1} + Y)}} - e^{-\gamma_{t}e^{-(\alpha_{n}Z_{n-1} + G_n  + \delta_{n})}}\big)\Big]}\\
  \shoveleft{   = \esp\Big[e^{-Z_{n-1}}\big(e^{-\gamma_{t}e^{-(\alpha_{n}Z_{n-1} + Y)}} - e^{-\gamma_{t}e^{-(Z_{n-1} + Y)}}\big)\ind_{\lb Z_{n-1} \leq 0 \rb}\Big]}\\
  + \esp\Big[e^{-Z_{n-1}}\big( e^{-\gamma_{t}e^{-(Z_{n-1} + Y)}} - e^{-\gamma_{t}e^{-(\alpha_{n}Z_{n-1} + Y)}}\big)\ind_{\lb Z_{n-1} > 0 \rb}\Big]\\
  \shoveright{+ \esp\Big[e^{-Z_{n-1}}\big( e^{-\gamma_{t}e^{-(\alpha_{n}Z_{n-1} + Y)}} - e^{-\gamma_{t}e^{-(\alpha_{n}Z_{n-1} + G_n  + \delta_{n})}}\big)\Big]}\\
  \shoveleft{\leq \esp\Big[e^{-Z_{n-1}}\big( e^{-\gamma_{t}e^{-(\alpha_{n}Z_{n-1} + Y)}} - e^{-\gamma_{t}e^{-(Z_{n-1} + Y)}}\big)\ind_{\lb Z_{n-1} \leq 0 \rb}\Big]}\\
  + \esp\Big[e^{-Z_{n-1}}\big( e^{-\gamma_{t}e^{-(Z_{n-1} + Y)}} - e^{-\gamma_{t}e^{-(\alpha_{n}Z_{n-1} + Y)}}\big)\ind_{\lb Z_{n-1} > 0 \rb}\Big]\\
  + \esp\Big[e^{-Z_{n-1}}\big( e^{-\gamma_{t}e^{-(\alpha_{n}Z_{n-1} + Y)}} - e^{-\gamma_{t}e^{-(\alpha_{n}Z_{n-1} + \frac{Y}{\varepsilon_{n}} + \delta_{n})}}\big)\Big].
\end{multline}
Inequality \eqref{cpn_Y_G_n} yielded the last line. The two first terms of that line become easy to estimate because for $\alpha \in \lb \alpha_{n}, 1\rb$:
$$
e^{-x}\esp\big[e^{-\gamma_{t}e^{-(\alpha x + Y)}}\big] = \frac{1}{\gamma_{t}}e^{-(1-\alpha)x}\big(1-e^{-\gamma_{t}e^{-\alpha x}}\big),\ x \in \R
$$
so that:
\begin{multline*}
  \esp\big[e^{-x}\big( e^{-\gamma_{t}e^{-(\alpha_{n}x + Y)}} - e^{-\gamma_{t}e^{-(x + Y)}}\big)\big] \\
  = \frac{1}{\gamma_{t}}\big((e^{-(1-\alpha_{n})x} - 1) + (e^{-\gamma_{t}e^{-x}} - e^{-(1-\alpha_{n})x}e^{-\gamma_{t}e^{-\alpha_{n}x}})\big)                  \\
  = \frac{1}{\gamma_{t}}\big((e^{-(1-\alpha_{n})x} - 1)(1 - e^{-\gamma_{t}e^{-\alpha_{n}x}}) + (e^{-\gamma_{t}e^{-x}} - e^{-\gamma_{t}e^{-\alpha_{n}x}})\big).
\end{multline*}
Now we can bound the first two terms of \eqref{eq_coupon_core:12} by putting the indicator functions together again:
\begin{multline*}
  \esp\Big[e^{-Z_{n-1}}\big\vert e^{-\gamma_{t}e^{-(Z_{n-1} + Y)}} - e^{-\gamma_{t}e^{-(\alpha_{n}Z_{n-1} + Y)}} \big\vert\Big] \\                                                                                                         		\leq \frac{1}{\gamma_{t}}\esp\Big[\big\vert e^{-(1-\alpha_{n})Z_{n-1}} - 1\vert + \vert e^{-\gamma_{t}e^{-Z_{n-1}}} - e^{-\gamma_{t}e^{-\alpha_{n}Z_{n-1}}} \big\vert\Big].
\end{multline*}
Note  first that  the function
$z \mapsto e^{-Z_{n-1}z}$ is Lipschitz on ${[0, 1-\alpha_{n}]}$, with Lipschitz
constant less than $\vert Z_{n-1} \vert (1 + e^{-(1-\alpha_{n})Z_{n-1}})$ almost
surely. Second, remark that the function $z \mapsto e^{-\gamma_{t}e^{-z}}$ is
$e^{-1}$ Lipschitz. Thanks to Cauchy-Schwartz inequality, we know that
$$
\underset{n \in \mathbb{N}^{*}}{\sup}\esp\big[\vert Z_{n-1} \vert e^{-\lambda Z_{n-1}}\big] < +\infty,\ \lambda \geq 0.
$$
Therefore, we can write:
\begin{align*}
  \esp\Big[e^{-Z_{n-1}}\big\vert e^{-\gamma_{t}e^{-(Z_{n-1} + Y)}} - e^{-\gamma_{t}e^{-(\alpha_{n}Z_{n-1} + Y)}}\big\vert\Big] \leq \frac{C}{\gamma_{t}}\frac{1}{n}\cdotp
\end{align*}
As to the rightmost term  of~\eqref{eq_coupon_core:12}, a change of variable gives:
\begin{align*}
  e^{-x}\esp\big[e^{-\gamma_{t}e^{-(\alpha_{n}x + \frac{Y}{\varepsilon_{n}} + \delta_{n})}}\big] = \varepsilon_{n}e^{-x}\int_{0}^{\infty}e^{-\varepsilon_{n}y}e^{-\gamma_{t}e^{-(\alpha_{n}x + y + \delta_{n})}}\dif y,\ x \in \R.
\end{align*}
We conclude with Lemma~\ref{lem_coupon_core:6}.
\end{proof}

By combining Lemmas \ref{lemma_A12}, \ref{lem_coupon_core:3}, \ref{lemma_C2}, and \ref{lem_coupon_core:5}, we have bounded each of the terms $A_{1}$, $A_{2}$ and $A_{3}$ appearing in our initial decomposition \eqref{cpn_dcptn}, thus proving theorem \ref{cpn_aroundI}.

  \subsection{Proof of theorem \ref{cpn_around0} - Bounding the error on $[0,1]$}

  We now have to find an upper-bound of $\esp[\Lo\Po_{t}h(Z_{n})]$ uniform with respect to~$t$. We rely once more on decomposition \eqref{cpn_dcptn}. First, we have for $A_{1}$:
  \begin{align*}\label{A1_bis}
	\vert A_{1} \vert = (1-K_{n})\vert\esp[g_{n}(Z_{n-1})]\vert & \leq \esp\Big[e^{-\beta_{n}Z_{n-1}}e^{-\gamma_{t}e^{-(\alpha_{n}Z_{n-1} + G_n  + \delta_{n})}}\Big] \\
											& \leq (1-K_{n})\esp[e^{-\beta_{n}Z_{n-1}}]                                                                      \\
											& \leq (1-K_{n})\numberthis,
  \end{align*}
  and thanks to \eqref{K_n}, we already know that this last term goes to $0$ as fast as $n^{-1}\log n$. Term $A_{2}$ is dealt with in the next lemma:

	\begin{lemma}\label{lemma_A2bis}
	  Let $g$ be defined as in \eqref{def_g}, with
	  $h\in \lipbis$. Then there exists $C > 0$ such that:
	  \begin{align*}
		\vert A_{2} \vert = \vert\esp[g(Z_{n})] - \esp[g(Z_{n-1})]\vert & \leq C\ \frac{\log n}{n},\ t \geq 0, n\geq 2.
	  \end{align*}
	\end{lemma}

	\begin{proof}
	  Recall that $Y$ is a random variable with exponential distribution
	  $\mathcal{E}(1)$ and independent of everything. By using that 
	  $$
	  z \mapsto e^{-\gamma_{t}e^{-z}}h'(z - t)
	  $$
	  is $(1 + e^{-1})$-Lipschitz, and bounded by $1$, we find that:
	  \begin{multline*}
		\vert\esp[g(Z_{n})] - \esp[g(Z_{n-1})]\vert                                                                                                                                      \\
		\begin{aligned}
		  &\leq \esp\Big[e^{-Z_{n}}\big\vert e^{-\gamma_{t}e^{-(Z_{n}+Y)}}h'(Z_{n} + Y - t) - e^{-\gamma_{t}e^{-(Z_{n-1}+Y)}}h'(Z_{n-1} + Y - t)\big\vert\Big] \\
		  & + \esp\big[\vert e^{-Z_{n}} - e^{Z_{n-1}}\vert e^{-\gamma_{t}e^{-(Z_{n-1} + Y)}}\vert h'(Z_{n-1} + Y - t) \vert\big]                          \\
		  &\leq (1 + e^{-1})\esp\big[e^{-Z_{n}}\vert Z_{n} - Z_{n-1}\vert\big] + (1 + e^{-1})\esp\big[\vert e^{-Z_{n}} - e^{-Z_{n-1}}\vert\big].
		\end{aligned}
	  \end{multline*}
	  To exploit the fact that $\tau_{i+1}^{n} \geq \tau_{i}^{n-1}$, recall that
	  $\tau_{1}^{n} = 1$ a.s. and introduce:
		$$
		Z_{n}' \coloneq \frac{1}{n-1}\sum_{i=1}^{n}\tau_{i}^{n} - \log(n-1) = \frac{1}{n-1}\sum_{i=1}^{n-1}\tau_{i+1}^{n} + \frac{1}{n-1} - \log(n-1),\ n \geq 2.
		$$
		Clearly $Z_{n}' \geq Z_{n-1}$ and $Z_{n}' \geq Z_{n}$ since
		$$
		Z_{n} - Z_{n}' = -\frac{1}{n(n-1)}\sum_{i=1}^{n}\tau_{i}^{n} + \log\Big(1 - \frac{1}{n}\Big) \leq 0.
		$$
		Therefore, by bounding $1 + e^{-1}$ by $2$:
    \begin{multline*}
      \frac{1}{2}\vert\esp[g(Z_{n})] - \esp[g(Z_{n-1})]\vert \\ \leq \esp\big[e^{-Z_{n}}\big\vert Z_{n} - Z_{n-1}\vert\big] + \esp\big[\vert e^{-Z_{n}} - e^{-Z_{n-1}}\big\vert\big] \\
{       \leq \esp\big[e^{-Z_{n}}(Z_{n}' - Z_{n})\big] + \esp\big[e^{-Z_{n}}(Z_{n}' - Z_{n-1})\big]}\\
     \shoveright{ + \esp\big[e^{-Z_{n}} - e^{-Z_{n}'}\big] + \esp\big[e^{-Z_{n-1}} - e^{-Z_{n}'}\big] }\\
        = \sum_{i=1}^{4}D_{i}.
    \end{multline*}
		We introduce the notation:
		$$
		Z_{n\bs i} \coloneq \frac{1}{n}\underset{\substack{j = 1 \\ j\neq i}}{\sum^{n}}\tau_{j}^{n} - \log n,\ i \in I_{n}.
		$$
		Bounding $-Z_{n}$ by $-Z_{n\bs i}$, we can write:
		\begin{align*}
		  D_{1} & = \frac{1}{n(n-1)}\sum_{i=1}^{n}\esp[e^{-Z_{n}}\tau_{i}^{n}] - \log\big(1 - \frac{1}{n}\big)\esp[e^{-Z_{n}}]        \\
												   & \leq \frac{1}{n(n-1)}\sum_{i=1}^{n}\esp[e^{-Z_{n}}\tau_{i}^{n}] + \frac{1}{n-1}                                     \\
												   & = \frac{1}{n(n-1)}\sum_{i=1}^{n}\esp[e^{-Z_{n\bs i}}]\esp[e^{-\frac{1}{n}\tau_{i}^{n}}\tau_{i}^{n}] + \frac{1}{n-1} \\
												   & \leq \frac{1}{n(n-1)}\sum_{i=1}^{n}\esp[e^{-Z_{n\bs i}}]\esp[\tau_{i}^{n}] + \frac{1}{n-1}                          \\
												   & = \frac{1}{n-1}\sum_{i=1}^{n}\esp[e^{-Z_{n\bs i}}]\frac{1}{n-i+1} + \frac{1}{n-1}\cdotp
		\end{align*}
		That last expectation is bounded with respect to $n$:
		\begin{align*}
		  \esp[e^{-Z_{n\bs i}}] = n\underset{\substack{j = 1                    \\ j\neq i}}{\prod^{n}}\frac{n-j+1}{ne^{1/n} - j + 1} & = \frac{nn!}{n-i+1}\underset{\substack{j = 1    \\ j\neq i}}{\prod^{n}}\frac{1}{ne^{1/n} - j + 1}\\
								& \leq \frac{nn!}{n-i+1}\underset{\substack{j = 1 \\ j\neq i}}{\prod^{n}}\frac{1}{n - j + 2}\\
								& = \frac{nn!}{n-i+1}\times\frac{n-i+2}{(n+1)!}   \\
								& = \frac{n-i+2}{n-i+1}\times\frac{n}{n+1}        \\
								& \leq 2.
		\end{align*}
		Consequently we have
		\begin{align*}
		  D_{1} &\leq \frac{1}{n-1} + \frac{1}{n(n-1)}\sum_{i=1}^{n}\esp[e^{-Z_{n\bs i}}]\esp[e^{-\frac{1}{n}\tau_{i}^{n}}\tau_{i}^{n}]\\
		&\leq \frac{1}{n-1} + \frac{2}{n-1}H_{n}\cdotp
		\end{align*}
		Those arguments also  apply  to $D_{2}$ and give almost the same bound:
		\begin{align*}
		 D_{2} & = \frac{1}{n-1} + \frac{1}{n-1}\sum_{i=1}^{n-1}\esp\big[e^{-Z_{n}}(\tau_{i+1}^{n} - \tau_{i}^{n-1})\big]                               \\
													 & = \frac{1}{n-1} + \frac{1}{n-1}\sum_{i=1}^{n-1}\esp[e^{-Z_{n\bs i}}]\esp\big[e^{-\tau_{i+1}^{n}}(\tau_{i+1}^{n} - \tau_{i}^{n-1})\big] \\
													 & \leq \frac{1}{n-1} + \frac{2}{n-1}\sum_{i=1}^{n-1}\frac{1}{n-i}                                                                        \\
													 & = \frac{1}{n-1} + \frac{2}{n-1}H_{n-1}\cdotp
		\end{align*}
    Use identity~\eqref{cpn_Z_n_asym} of lemma \ref{cpn_exp_bnd} to derive the asymptotic behavior of $D_{3}$ and $D_{4}$ and conclude the proof.
	  \end{proof}

	  As for $A_{3}$, we use again the fact that
	  $z\mapsto e^{-Z_{n-1}z}$ is
	  ${\vert Z_{n-1}\vert (e^{-\beta_{n}Z_{n-1}} + e^{-Z_{n-1}})}$-Lipschitz
    continuous almost surely when $z \in [\beta_{n}, 1]$:
	  \begin{multline*}
		\Big\vert \esp\Big[(e^{-\beta_{n}Z_{n-1}} - e^{-Z_{n-1}})e^{-\gamma_{t}e^{-(\alpha_{n}Z_{n-1} + G_n  + \delta_{n})}}  h'(Z_{n-1} + Y - t)\Big]\Big\vert                                                         \\
		\begin{aligned}
		  & \leq \esp\big[\vert e^{-\beta_{n}Z_{n-1}} - e^{-Z_{n-1}}\vert\big]                            \\
		  & \leq (1-\beta_{n})\esp\Big[\vert Z_{n-1}\vert (e^{-\beta_{n}Z_{n-1}} + e^{-Z_{n-1}})\Big] \\
		  & \leq C(1-\beta_{n}).
		\end{aligned}
	  \end{multline*}
	  We combine the previous inequality, Lemma \ref{lemma_A2bis} as well as inequality \ref{A1_bis} to get the bound of theorem \ref{cpn_around0}.


\section*{Declarations}

\begin{itemize}
\item Funding : LD acknowledges support from the French National Agency for Research (ANR) via the project n°ANR-22-PEFT-0010 of the France 2030 program PEPR réseaux du Futur.

\end{itemize}

\end{document}